\documentclass[12pt]{article}

\newtheorem{thm}{Theorem}[section]
\newtheorem{lemma}{Lemma}[section]
\newtheorem{cor}{Corollary}[section]

\usepackage{amssymb}
\usepackage{amsfonts}
\usepackage{latexsym}

\renewcommand{\phi}{\varphi}
\renewcommand{\L}{\mathcal L}
\renewcommand{\d}{\delta_{\{M_n\}} (s)}
\renewcommand{\P}{\mathcal P}

\newcommand{\R}{\mathbb R}
\newcommand{\C}{\mathbb C}
\newcommand{\F}{\mathcal F}
\newcommand{\M}{\mathcal M}

\newcommand{\Cl}{{\rm Clos}_{C^0_T} \P}
\newcommand{\rpW}{\rho_{p, W}}

\newcommand{\lW}{\lambda_{W}}

\newcommand{\intR}{\int_{\R}}

\title{Asymptotics of the Fourier and Laplace transforms 
in weighted spaces of analytic functions
\thanks
{Supported in part by the Israel Science Foundation of
the Israel Academy
of Sciences and Humanities under Grants Nos 93/97-1,
37/00-1.}}

\author{V. Matsaev, M. Sodin}

\date{}

\begin{document}
\maketitle

\centerline{\em In memory of Evsei Dyn'kin (1949--1999)}

\begin{abstract}
We study the asymptotics near the origin of the Fourier transform in 
weighted Hardy spaces of analytic functions in the upper half-plane, and
of the Laplace transform in weighted spaces of entire functions of zero
exponential type. 

These results are applied to two closely related problems posed by
E.~Dyn'kin: we find the asymptotics of the depth of zero in
non-quasianalytic Denjoy-Carleman classes, and of the
exact Levinson-Sj\"oberg majorant.  

\end{abstract}

\bigskip

\section{Main results}
\setcounter{equation}{0}

Let $H^p$, $1\le p \le \infty$, be the Hardy spaces of analytic functions
in the upper half-plane $\C_+$, and let $W$ be a non-vanishing analytic
function in $\C_+$. The {\em weighted Hardy spaces} $H^p(W)$ are defined
as
follows:
$$
H^p(W) = \left\{
f:\, f\ {\rm is\ analytic\ in}\ \C_+\,,
\quad Wf \in H^p 
\right\},
$$
$$
||f||_{H^p(W)} \stackrel{def}= ||Wf||_{H^p}\,.
$$
Throughout this paper, we assume that the weight function $W$ is an outer 
function, that is
\begin{equation}
W(z) = \exp \left\{ \frac{1}{\pi i} \intR 
\left[ \frac{1}{t-z} - \frac{t}{t^2+1} 
  \right] \phi (t)\, dt\right\}.
\label{1.1}
\end{equation}
About the {\em logarithmic weight} $\phi$ we always assume that
{\em \smallskip\par\noindent (i) $\phi$ is non-negative and even on $\R$;

\smallskip\par\noindent (ii) $\phi$ increases on $(0, \infty)$ and 
$$
\lim_{t\to\infty} \frac{\phi (t)}{\log t} = \infty\,;
$$

\smallskip\par\noindent (iii) 
$$
\intR \frac{\phi (t)}{t^2+1}\, dt < \infty\,;
$$

\smallskip\par\noindent (iv) at least one of the following two conditions
holds true:

\smallskip\par\noindent \qquad (iv-a) the function $\tau \mapsto
\phi(e^\tau)$
is
convex;

\smallskip\par\noindent \qquad (iv-b) the function $\phi (t)$ is concave
on
$(0,\infty)$ and 
$$
\lim_{t\to \infty} t\phi'(t+0) = \infty\,.
$$
}
(In the case (iv-a), the latter limit is also infinite
due to condition (ii)).

Below, we tacitly assume that $\phi$ is a $C^2$-function, and $\phi (0) =
0$; however these assumptions can be easily dropped. 

Functions from the space
$H^p(W)$ rapidly decrease along the real axis, and we can define the
(inverse) Fourier transform
$$
(\F^{-1} f)(s) = \frac{1}{\sqrt{2\pi}} \intR f(x) e^{-isx}\, dx\,. 
$$
Our assumptions guarantee that $\F^{-1} f$ is a $C^\infty$-function
vanishing for $s\le 0$ (the integration can be shifted to the line $\R+iy$
for an arbitrarily large positive $y$). The question we are interested in
is:
\begin{itemize}
\item[$\bullet$]
{\em how deep is the zero of the function $(\F^{-1}f)(s)$ at $s=0$?} 
\end{itemize}
We 
measure its depth by
$$
\rpW (s) \stackrel{def}= \sup_{||f||_{H^p(W)}\le 1}\,
\left|(\F^{-1}f)(s)\right|\,,
$$
and we compare this quantity with the (upper) Legendre transform
\begin{equation}
Q(s) \stackrel{def}= \sup_{y>0}
\left[ \log |W(iy)| - ys 
\right] 
= \sup_{y>0} \left[ \frac{2y}{\pi} \int_0^\infty 
\frac{\phi(t)dt}{t^2+y^2} - ys 
\right].
\label{1.2}
\end{equation}
After a change of variables, one can express this function in another
form. Let $q(y) = \log |W(iy)|$. In our assumptions (i)--(iv), the
function $q(y)$ is a concave increasing function, and $q'(y)$ steadily
tends to zero when $y$ goes to infinity (see Corollary~3.2 below). Then
\begin{equation}
Q(q'(y)) = q(y) - yq'(y) = \frac{4y^3}{\pi} 
\int_0^\infty \frac{\phi (t)\,dt}{(t^2+y^2)^2}\,.
\label{temp*}
\end{equation}
Loosely speaking, the function $Q(q'(y))$ has the same growth as $\phi
(y)$; more precisely, once can check using (\ref{temp*}) that  
$$
\frac{4}{3\pi} \phi (y) \le Q(q'(y)) \le
\frac{16}{3\pi} y^3 \int_y^\infty \frac{\phi (t)}{t^4}\, dt\,.
$$

\begin{thm}
Let the logarithmic weight $\phi $ satisfy conditions (i)-(iv), and let
\smallskip\par\noindent (v)
$$
\lim_{t\to\infty} \frac{t\phi'(t)}{ 
\left(t^3 \int_t^\infty \frac{\phi (\xi)}{\xi^4}\, d\xi \right)^{2/3}
} = +\infty\,.
$$  
Then for $s\to 0$
\begin{equation}
\log \rpW (s) = - (1+o(1)) Q(s)\,, 
\label{*}
\end{equation}
and 
\begin{equation}
\left( \F^{-1}\frac{1}{W}\right)(s) = (1+o(1))\sqrt{Q''(s)}\, e^{-Q(s)}\,.
\label{**}
\end{equation}
Furthermore, if additionally, $\phi$ has a positive lower order, i.e 
\par\noindent (vi) 
$$
\liminf_{t\to\infty} \frac{\log \phi (t)}{\log t} > 0\,,
$$
then for $s\to 0$
\begin{equation}
\log \rpW (s) = -Q(s) + O(\log Q(s))\,.
\label{***}
\end{equation}
\end{thm}

\medskip Condition (v) is not very restrictive since generally
speaking $t\phi '(t)$ has the same order of growth as the expression 
$$t^3\int_t^\infty \frac{\phi (s)}{s^4}ds\,.$$ For example, if $\phi $ is
concave (or, less restrictively, the function $\phi '(t)/t$ does not
increase), and if the function $\phi (t)\log^{-3}t$ steadily
increases with $t$ to $+\infty$, then by a straightforward
inspection
\footnote{
Indeed, 
$$
t^3 \int_t^\infty \frac{\phi (\xi)}{\xi^4}\,d\xi 
= \frac{1}{3} \left\{ \phi (t) + t^3 \int_t^\infty
\frac{\phi'(\xi)}{\xi^3}\,d\xi
\right\} \le \frac{1}{3} \big\{
\phi (t) + t\phi'(t) \big\}\,.
$$
Since the function $\phi (t) \log^{-3}t$ increases, we have  
$$
t\phi '(t) \ge \frac{3\phi (t)}{\log t}\,,
$$
and 
$$
\frac{t\phi'(t)}{\phi^{2/3}(t)} \ge 3 \frac{\phi^{1/3}(t)}{\log t}
\uparrow \infty\,.
$$
This yields condition (v)}
condition (v) holds. We do not know whether condition (v) is
really needed for asymptotics (\ref{*}) and (\ref{***}). 

In a special case, when the weight function is an entire function of genus
zero with zeros on the negative imaginary semi-axes:  
$$
W(z) = \prod_k \left( 1 + \frac{z}{ia_k}\right)\,,
\qquad \sum_{k} \frac{1}{a_k} <\infty\,,
$$
asymptotics (\ref{*}), (\ref{**}) and (\ref{***}) hold without any
additional restrictions. This can be deduced from a result of Hirschman
and Widder \cite[Chapter~V,~\S3]{HW}.

\medskip A dual problem to the one considered above is:
\begin{itemize}
\item[$\bullet$]
{\em to estimate the
growth of the Laplace transform 
$$
(\L f)(\zeta) = \int_0^\infty f(x)e^{-\zeta x}\, dx\,,
$$ 
of an entire function $f$ of zero exponential type, when $\zeta$
approaches the origin.}
\end{itemize}
Let $W$ be an outer function (\ref{1.1}). Define
the space $B(W)$ of entire functions: 
$$
B(W) = \left\{f:\, f\ {\rm is\ entire}, \quad
f/W, f^*/W \in H^\infty
\right\},
$$
where $f^*(z)=\overline{f(\bar z)}$,
$$
||f||_{B(W)} \stackrel{def}= \sup_{t\in \R} \frac{|f(t)|}{|W(t)|}\,.
$$
Since the function $W$ is outer, $B(W)$ consists of entire functions of
zero exponential type. The Laplace transform $(\L f)(\zeta)$ is analytic
in the right half-plane ${\rm Re}\, \zeta>0$, and since $f$ has zero
exponential type, we can turn the integration line in any direction and
thus obtain that $\L f$ is analytic in ${\overline \C} \setminus \{0\}$
and vanishes at infinity. Set
$$
\lW (s) \stackrel{def}=
\sup_{||f||_{B(W)}\le 1} \, \max_{|\zeta|=s} |(\L f)(s)|\,.
$$

\begin{thm}
Let the logarithmic weight $\phi$ satisfy conditions (i)-(iii), (iv-a) and
(iv-b). Then for $s\to 0$
\begin{equation}
\log \lW (s) = (1+o(1))Q(s)\,.
\label{lambda1}
\end{equation}
If additionally condition (vi) holds, then 
\begin{equation}
\log \lW (s) = Q(s) + O(\log Q(s))\,.
\label{lambda2}
\end{equation}
\end{thm}

\medskip We do not know whether condition (iv-a) is needed for asymptotics 
(\ref{lambda1}).

\medskip It is possible to introduce the scale of de Branges-type spaces
$B^p(W)$, $1\le p \le \infty$. In all these spaces the asymptotic
relations (\ref{lambda1}) and (\ref{lambda2}) continue to hold.

\medskip The proofs are based on the classical Laplace method of
asymptotic evaluation of integrals. Let us mention briefly another
situation when the technique developed here could be useful. 
Let us drop condition (iii) and consider the so-called {\em quasianalytic
weights} $\phi (t)$ such that 
$$
\intR \frac{\phi (t)}{t^2+1}\, dt = \infty\,,
$$
but
$$
\intR \frac{\phi (t)}{1+|t|^3}\, dt < \infty\,.
$$
Then one can define an analytic weight function $W(z)$ using the modified
Schwartz integral:
$$
W(z) = \exp 
\left\{ \frac{1}{\pi i} \intR \frac{(tz+1)^2 \phi (t)}{(t^2+1)^2 (t-z)}\,
dt \right\}\,,
$$ 
and look at the asymptotics of the Fourier transform
$$
\left( \F \frac{1}{W} \right)(s) = 
\frac{1}{\sqrt{2\pi}} \intR \frac{e^{isx}}{W(x)}\, dx\,,
\qquad s\to \infty\,,
$$
and of the Laplace transform
$$
\left( \L \frac{1}{W} \right)(s) = 
\int_0^\infty \frac{e^{sy}}{W(iy)} \, dy\,,
$$
when $s\to\infty$ in $\mathbb C$. These functions are closely related to
the study of primary ideals in weighted convolution algebras on $\mathbb
R$ with quasianalytic weights (cf. Borichev's work \cite{Borichev2} and
references therein). 

\bigskip
The authors thank Alexander Borichev and Iossif Ostrovskii for 
spotting a number of blunders and for numerous helpful remarks and 
suggestions.

\section{Some applications}
\setcounter{equation}{0}

\subsection{Depth of zero in non-quasianalytic classes} Let
$C\{M_n\}$ be a
non-quasianalytic Denjoy-Carleman class of $C^\infty$-functions on $\R$
such that
$$
\sup_{\R} |g^{(n)}| \le M_n\,,
\qquad n=0, 1, 2, \, ...\, .
$$
We are interested in 
\begin{itemize}
\item[$\bullet$]
{\em 
the asymptotic behaviour of the quantity
$$
\d = \sup \left\{ |g(s)|:\,
g\in C\{M_n\}, \quad g^{(n)}(0)=0\,, \ n=0, 1, 2, \,...
\right\}\,,
$$
when $s\to 0$.}
\end{itemize}

This problem was raised and treated by Th.~Bang in the paper \cite{Bang}; 
a related problem in quasianalytic classes was considered earlier by 
Carleman \cite[pp. 24--27]{Carleman}. Bang observed
that if the sequence $\{M_n\}$ grows sufficiently slowly (but still
non-quasianalytically), then the estimate obtained from the Taylor formula
\begin{equation}
\d \le \inf_{n\ge 0} \frac{M_ns^n}{n!}
\label{2.1}
\end{equation}
leads to a very crude bound which does not capture the transition from
non-quasianalyticity to quasianalyticity. For example, if 
\begin{equation}
M_n=n! (\log n)^{n(1+\beta)}\,, \qquad \beta>0\,,
\label{2.2}
\end{equation}
then (\ref{2.1}) yields 
$$
\log \frac{1}{\d} \ge e^{c s^{-1/(\beta +1)}}\,.
$$ 
This bound does not blow up when $\beta\to 0$ as it should. The
bound obtained by Bang by an ingenious and elementary method is
$$
\log \frac{1}{\d} \ge e^{c s^{-1/\beta }}\,.
$$ 

The problem of finding the asymptotic behaviour of $\d$ was raised anew by
Dyn'kin \cite{D2} in the seventies at a time when he was not aware of
Bang's work. A lower bound for $\d$ follows from
the construction of V.~P.~Gurarie \cite{Gurarie} and Beurling 
\cite{Beurling}. The upper
bound for $\d$ was obtained by Volberg \cite{Volberg} and 
Dyn'kin \cite{D3} by different methods. Borichev has a related result
in \cite{Borichev}. These results provide good bounds for
$|\log\d |$ but are far from being the asymptotics. For example, in the
case (\ref{2.2}), the best they could achieve is
$$
c s^{-1/\beta} \le \log |\log \d | \le Cs^{-1/\beta}\,,
\qquad s\le s_0\,.
$$
The discrepancy between the upper and lower bounds worsens if the 
class $C\{M_n\}$ is ``less'' non-quasianalytic. 

Here, using Theorem~1.1 we  obtain a much sharper result. Let
$$
T(r) = \sup_{n\ge 0} \frac{r^n}{M_n}
$$
be A. Ostrowski's function. Then the function 
\begin{equation}
\phi (t) = \log T(|t|) = \sup_{n\ge 0} 
\left[ n\log |t| - \log M_n \right]
\label{+}
\end{equation}
automatically satisfies conditions (i)-(iii) and (iv-a). 

\begin{lemma} Let $W$ be the weight function defined by (\ref{1.1})
with $\phi $ given by (\ref{+}). Then
\begin{equation}
\sqrt {2\pi}\, \rho_{1, W} (s) \le \d \le \frac{e}{\sqrt{2\pi}}\, 
s \rho_{\infty, W}(s)\,.
\label{2.3}
\end{equation}
\end{lemma}

\smallskip\par\noindent{\em Proof of Lemma~2.1:} 
We repeat Carleman's argument \cite{Carleman} with minor variations. Let
$g\in C\{M_n\}$, and $g^{(n)}(0)=0$, $n\in \mathbb Z_+$.
Without loss of generality, we assume that $g(s)\equiv 0$ for $s\le 0$, 
and consider the Fourier transform
$$
(\F g)(z) = \frac{1}{\sqrt {2\pi}}
\int_0^\infty g(s) e^{isz}\, ds
$$
analytic in the upper half-plane. Then
$$
(\F g)(z) = \frac{(-1)^n}{(iz)^n \sqrt{2\pi}} 
\int_0^\infty g^{(n)}(s) e^{isz}\, ds\,,
$$
whence
$$
|(\F g)(z)| \le \frac{1}{\sqrt{2\pi} yT(|z|)}\,,
\qquad z\in\C_+\,,
$$
or
$$
||(\F g)(z+i\tau)||_{H^\infty(W)} \le \frac{1}{\sqrt{2\pi}\, \tau}\,,
\qquad \tau>0\,,
$$
where the weight $W$ is defined by (\ref{+}). Then
$$
g(s) = \frac{1}{\sqrt{2\pi}} \intR (\F g)(x+i\tau) e^{-is(x+i\tau)}\,
dx\,,
$$
and
$$
|g(s)| \le \frac{e^{s\tau}}{\sqrt{2\pi}\, \tau} \rho_{\infty, W}(s)\,.
$$
Choosing here the optimal value $\tau=1/s$, we obtain the second half of
(\ref{2.3}).

Now, let $f\in H^1(W)$, and $||f||_{H^1(W)}\le 1$. Then $(\F^{-1} f)(s)$
vanishes for $s\le 0$, and
$$
\left| (\F^{-1}f)^{(n)} (s)\right|
\le \frac{1}{\sqrt{2\pi}} \intR \frac{|x|^n}{T(|x|)}\, |f(x)| T(|x|)\, dx
\le \frac{M_n}{\sqrt{2\pi}}\,.
$$
Therefore, $\sqrt{2\pi} \F^{-1}f \in C\{M_n\}$ and the first
inequality in (\ref{2.3}) follows. $\Box$

\medskip
Combining this lemma with Theorem~1.1, we get

\begin{thm}
If the logarithm of Ostrowski's function $\phi(r) = \log T(r)$
satisfies condition (v), then
$$
\log \d  = -(1+o(1)) Q(s)\,, \qquad 
s\to 0\,,
$$
where $Q(s)$ is the upper Legendre transform (\ref{1.2}).

If $\phi (r)$ also satisfies condition (vi), then a sharper asymptotic
relation 
$$
\log \d  = -Q(s) + O(\log Q(s))\,, \qquad 
s\to 0\,,
$$
holds.
\end{thm}

\subsection{Distance from the polynomials to an
imaginary exponent in weighted spaces. } Our second application pertains 
to the weighted polynomial approximation. 

Let $T(r)=e^{\phi (r)}$, where $\phi$ satisfies assumptions (i)--(iv), and 
let $C_T^0$ be the weighted 
space of continuous functions $h$ on $\R$ such that
$$
\lim_{t\to \infty} \frac{h(t)}{T(|t|)} = 0\,,
$$
and
$$
||h||_{C^0_T} = \sup_{t\in \R}\, \frac{|h(t)|}{T(|t|)}\,.
$$
Let $X=\Cl$ be the closure of algebraic polynomials $\P$ in $C_T^0$.

Due to condition (iii), the polynomials are not dense in $C_T^0$ and 
\begin{equation}
\Cl \subset \mathcal E_0 \cap C_T^0\,,
\label{2.4}
\end{equation}
where $\mathcal E_0$ is a space of all entire functions of zero 
exponential type (cf. \cite[Section~VI]{Koosis}). It has been 
known for a long time that 
under condition (iv-a) plus some mild regularity of $\phi$ there is 
always an equality sign in (\ref{2.4}) (Khachatrian, Koosis,
Levinson-McKean). 
Recently, Borichev proved in \cite{Borichev4} that this always holds  
whenever $\phi$ satisfies conditions (i)-(iii) and (iv-a).
We shall not use this result.

Let $e_s(t)=e^{ist}$. 
Here, we are interested in 
\begin{itemize}
\item[$\bullet$]
{\em 
the asymptotics of the quantity
$$
d_T(s) \stackrel{def}= {\rm dist}_{C_T^0} (X, e_s) = 
{\rm dist}_{C_T^0} (\P, e_s)\,,
$$
when $s\to 0$.}
\end{itemize}
A well-known duality argument links this question to 
the previous one:

\begin{lemma}
Let $W$ be an outer function 
(\ref{1.1}) with the boundary values $|W(t)|=T(|t|)$, $t\in\R$. Then
$$
\sqrt{2\pi}\, \rho_{1,W}(s) \le d_T(s) \le \frac{e}{\sqrt{2\pi}}\, s
\rho_{\infty, W} (s)\,.
$$
\end{lemma}

\smallskip{\em Proof of Lemma~2.2:}
We identify the dual space $\left( C_T^0\right)^*$ with the space of 
complex-valued measures $\mu$ on $\R$ such that 
$$
||\mu||_{T} \stackrel{def}= \intR T(|t|)\, d\mu (t) < \infty\,.
$$
Then by the Hahn-Banach theorem
\begin{equation}
d_T(s) = \sup_{\mu\in\P^\perp, \, ||\mu||_T\le 1} |\mu (e_s)|
= \sqrt{2\pi}\, \sup_{\mu\in\P^\perp, \, ||\mu||_T\le 1} |(\F\mu)(s)|\,,
\label{2.5}
\end{equation}
where
$$
(\F\mu)(s) = \frac{1}{\sqrt{2\pi}} \intR e^{ist}\, d\mu (t)
$$
is the Fourier transform of the measure $\mu$, and
$\P^\perp \subset \left( C_T^0\right)^*$ is the annihilator of the 
polynomials. Clearly, $(\F\mu)(s)$ is a $C^\infty$-function on $\R$ such 
that
$$
(\F\mu)^{(n)}(0) = \frac{i^n}{\sqrt{2\pi}}
\intR t^n\, d\mu (t) = 0\,,
\qquad n\in\mathbb Z_+\,,
$$
and
\begin{eqnarray*}
\left|(\F\mu)^{(n)}(s)\right|
&=& \left| \frac{1}{\sqrt{2\pi}} \intR (it)^n e^{ist}\, d\mu (t) 
\right| \\ \\
&\le& \frac{1}{\sqrt{2\pi}} \intR |t|^n \, |d\mu (t)| \\ \\
&= & \frac{1}{\sqrt{2\pi}} \intR \frac{|t|^n}{T(|t|)} T(|t|) \, |d\mu (t)| 
\\ \\
&\le& \frac{1}{\sqrt{2\pi}} \sup_{r\ge 0} \frac{r^n}{T(r)} 
=: \frac{1}{\sqrt{2\pi}} M_n
\end{eqnarray*}
(the function $T(r)$ is Ostrowski's function for the sequence $M_n$ 
defined in the last line). Juxtaposing this estimate with (\ref{2.5}) and
using the Lemma~2.1, we obtain
$$
d_T(s) \le \d \le \frac{e}{\sqrt{2\pi}} s\rho_{\infty, W}(s)\,.
$$

In order to get the lower bound for $d_T(s)$, we define the measure
$d\mu (t) = f(t)\,dt$, where $f\in H^1(W)$ and $||f||_{H^1(W)} \le 1$.
Then $||\mu||_T = ||f||_{H^1(W)} \le 1$. Since $\F^{-1} f$ always has a 
zero of infinite order at the origin, this measure annihilates the 
polynomials. Therefore, by (\ref{2.5})
$$
d_T(s) \ge \sup_{||f||_{H^1(W)}\le 1} \left| (\F^{-1}f)(s)\right|
= \sqrt{2\pi} \rho_{1, W} (s)\,,
$$
completing the proof. $\Box$

\smallskip Combining this lemma with Theorem~1.1, we obtain

\begin{thm}
Let the function $\phi = \log T$ satisfy condition (v).
Then
$$
\log d_T (s) = -(1+o(1))Q(s)\,,
\qquad s\to 0\,, 
$$
where $Q$ is the upper Legendre transform (\ref{1.2}). 
If $\phi $ also satisfies condition (vi), then 
$$
\log d_T (s) = -Q(s) + O(\log Q(s))\,,
\qquad s\to 0\,.
$$
\end{thm}

\smallskip The same asymptotics holds in the weighted 
$L^p$-spaces on $\R$ under the same assumptions on the weight $T$.

\subsection{Exact Levinson-Sj\"oberg majorant} 
Let $S=\{\zeta=\xi+i\eta: |\xi|<1, |\eta|<1\}$. 
The Carleman-Levinson-Sj\"oberg theorem 
says that {\em the set of analytic functions in $S$ such that
\begin{equation}
|F(\zeta)| \le \M (|\xi|)
\label{2.6}
\end{equation}
is locally uniformly bounded in $S$ provided that the majorant $\M(\xi)$ 
does not increase on $(0,1)$, and}
\begin{equation}
\int_0 \log \log \M(\xi)\, d\xi < \infty\,.
\label{2.7}
\end{equation}
In this form the theorem was proved independently by Levinson 
\cite{Levinson} and Sj\"oberg \cite{Sjoberg}. However a decade before 
Carleman proved an equivalent result in \cite {Carleman2}. 
It is easy to see that the result persists for analytic functions in the
punctured square $S^*=S\setminus \{0\}$ which satisfy conditions
(\ref{2.6}) and (\ref{2.7}); i.e. the corresponding family is locally
uniformly bounded in $S^*$.
 
In \cite{D2}, \cite{D3}, E.~Dyn'kin raised the question of 
\begin{itemize}
\item[$\bullet$]
{\em finding the growth of the exact majorant
$$
\M^* (s) = \sup_F \, \max_{|\zeta|=s} |F(\zeta)|\,, \qquad
s\to 0\,,
$$
where the supremum is taken over all analytic functions $F$ in the 
$S^*=S\setminus \{0\}$ satisfying conditions (\ref{2.6}), (\ref{2.7}).}
\end{itemize} 
An upper bound for $\M^*(s)$ can be obtained using Domar's approach
\cite{Domar_a}, \cite{Domar_b} (cf. \cite[Section~VII~D7]{Koosis}). 
A duality argument developed by the first-named author in 
\cite{Matsaev} shows that the Carleman-Levinson-Sj\"oberg theorem is {\em 
equivalent} to the Denjoy-Carleman quasianalyticity theorem. Later, this 
fact was re-discovered by Dyn'kin in \cite{D1}. This leads to a
lower bound for $\M^*(s)$ in terms of $\delta_{\{M_n\}}(s)$ with a
properly chosen non-quasianalytic sequence $\{M_n\}$.
The upper and lower bounds obtained in this way are not tight and there is
a gap between them. 
See the discussion and summary of known results in the survey papers 
\cite{D3} and \cite[pp. 69--71]{Nik}. 

Here, we shall not use
these results, but directly exploiting  Theorem~1.2 we find
an asymptotics for $\log \M^* (s)$ when $s\to 0$. 

Let
\begin{equation}
\phi (r) \stackrel{def}= \inf_{\xi>0}
\left[ \log \M(\xi) + r\xi \right]
\label{2.8}
\end{equation}
be the (lower) Legendre transform of $\log \M(\xi)$. We assume that the 
majorant $\M$ grows sufficiently fast: for each $N< \infty$
\begin{equation}
\lim_{\xi\to 0} \xi^N \M(\xi) = \infty\,.
\label{2.9}
\end{equation}
Then the logarithmic weight function $\phi (r) $ automatically satisfies
conditions (i)--(iii) and (iv-b): condition (iii) follows from
(\ref{2.7}), see \cite[Section~VIIID2]{Koosis}, and (\ref{2.9}) yields
condition (ii).

To ensure condition (iv-a) for the function $\phi$, we assume that
\begin{equation}
\mbox{the functions} \quad s\mapsto \log\M(e^{-s}) \quad \mbox{and} \quad 
\xi\mapsto
\log\M(\xi) \quad \mbox{are convex.}
\label{N} 
\end{equation}
Then the function $\phi(e^\tau)$ is convex. Indeed, set $m(\xi) =
\log\M(\xi)$ and assume without loss of generality that $m$ is a
$C^2$-function. Let $\xi=\xi_r$ be the unique solution of the equation
$m'(\xi)=-r$. Then
$$
r\phi''(r) + \phi'(r) =
\frac{-m'(\xi)}{-m''(\xi)} + \xi 
= \frac{m'(\xi) + \xi m''(\xi)}{m''(\xi)} \stackrel{(\ref{N})}> 0\,,
$$
which yields convexity of the function $\tau\mapsto \phi(e^\tau)$.

\begin{lemma}
Let the majorant $\M(\xi)$ satisfy conditions (\ref{2.7}) and
(\ref{2.9}). Then
$$
\lambda_{W^*} (s) \le \M^*(s) \le C\, \lambda_{W}(s)\,,
\qquad s\le 1/2\,,
$$
where $W$ is the outer function (\ref{1.1}) constructed by the logarithmic 
weight $\phi$ from (\ref{2.8}), $W^*(z)= W(z)/(z+i)^2$,
and $\lW$ is the same as in the
Theorem~1.2.
\end{lemma}

\smallskip\par\noindent{\em Proof of Lemma~2.3: }
Let $\gamma$ be a simple closed curve in $S^*$ which surrounds the origin, 
and let $\Omega_o$ and $\Omega_i$ be the outer and inner components of 
$\overline\C \setminus\gamma$, i.e. $\Omega_o$ contains infinity, and
$\Omega_i$ contains the origin. If $F$ is analytic in $S^*$, then
$$
\frac{1}{2\pi i} \int_\gamma \frac{F(w)w\, dw}{\zeta(w-\zeta)}
= \left\{
\begin{array}{ll}
F_o(\zeta), \quad &\zeta\in\Omega_o, \\ \\
F_i(\zeta), \quad &\zeta\in\Omega_i,
\end{array}
\right.
$$
where $F_o$ is analytic in $\overline\C \setminus\{0\}$ and decays at 
infinity at least as $O(|\zeta|^{-2})$, $F_i$ is analytic 
in $S^*$ and has at most a simple pole at the origin, and
$$
F(\zeta) = F_i(\zeta) - F_o(\zeta)\,,
\qquad \zeta\in S^*\,.
$$

By the Carleman-Levinson-Sj\"oberg theorem,
$$
\max_\gamma |F(w)| \le C_{\gamma, \M}\,,
$$
where $C_{\gamma, \M}$ depends on the function $\M$ and on ${\rm
dist}(\gamma, \partial S^*)$ but is independent of $F$. Therefore,
$$
|F_i(\zeta)| \le \frac{C}{|\zeta|}, \qquad |\zeta|\le \frac{1}{2}\,.
$$

Now, for $|\zeta|\le 1/2$
$$
|F_o(\zeta)| \le |F(\zeta)| + |F_i(\zeta)| \le \M(|\xi|) +
\frac{C}{|\zeta|} \le C\M(|\xi|)
$$
(the constants $C=C_\M$ may vary from line to line). For $|\zeta|\ge 1/2$,
we have
$$
|F_o(\zeta)| \le \frac{C}{|\zeta|^2}\,.
$$
So that,
$$
|F_o(\zeta)| \le \frac{C \M(|\xi|)}{1+\eta^2}\,,
\qquad \zeta\in\overline{\C}\setminus \{0\}
$$
(we put here $\M(\xi)=\M(1-0)$ for $\xi\ge 1$).

Since $F_o$ is analytic at $\overline{\C}\setminus\{0\}$ and vanishes at
infinity, it is a Laplace transform of an entire function $f$ of zero
exponential type:
$$
F_o(\zeta) = \int_0^\infty f(x)e^{-\zeta x}\, dx\,, 
\qquad \xi={\rm Re}\, \zeta >0\,.
$$ 
By the inversion formula
$$
f(x)e^{-\xi x} = \frac{1}{2\pi } \intR e^{i\eta x} F_o(\xi+i\eta)\,
d\eta\,,
\qquad x>0\,, \quad \xi>0\,,
$$
whence 
$$
|f(x)| \le C\, \inf_{\xi>0} \left[ \M(\xi)e^{\xi x}\right]
= Ce^{\phi (x)}\,,
\qquad x>0\,.
$$
The same  estimate holds for $x<0$. Therefore, $f\in B(W)$,  and for 
$s\le 1/2$
$$
\max_{|\zeta|=s} |F(\zeta)| 
\le \max_{|\zeta|=s} |F_o(\zeta)| + \max_{|\zeta|=s} |F_i(\zeta)|
\le C\left( \lambda_{\infty, W} (s) + \frac{1}{s} \right) 
\le C\lambda_{\infty, W}(s)\,, 
$$
proving the upper bound for $\M^*(s)$.

To get the lower bound, we take $f\in B(W^*)$. Then the Laplace
transform $\L f$ is analytic in $\overline{\C}\setminus\{0\}$, and for
$\xi>0$
\begin{eqnarray*}
\left| (\L f)(\xi+i\eta)\right| &\le& 
\int_0^\infty
e^{-\xi x + \phi (x)}\, \frac{|f(x)|(1+x^2)}{e^{\phi (x)}}\,
\frac{dx}{1+x^2} 
\\ \\
&\le& \pi\, ||f||_{B(W^*)} \exp\{\sup_{x>0}[\phi (x) - \xi x]\}
= \pi ||f||_{B(W^*)} \M(|\xi|)\,.
\end{eqnarray*}
The same bound holds for $\xi <0$. Therefore, $\M^*(s)\ge
\pi \lambda_{W^*}(s)$ completing the proof of the Lemma. $\Box$

\medskip 
Set $Q^*(s)=\sup_{y>0} [\log|W^*(iy)|-ys]$. It can be easily seen 
that $$Q(s)-o(Q(s)) \le Q^*(s) \le Q(s)\,.$$ 
Indeed, if the maximum in
(\ref{1.2}) is attained at a point $y_s$ (this point is unique and tends 
to $\infty$ as $s$ approaches the origin), then
\begin{eqnarray*}
Q(s) \ge Q^*(s) &\ge& \log|W(iy_s)| - y_s s - 2\log(1+y_s)
\\ \\
&=& Q(s) - 2\log(1+y_s) = Q(s) - o(Q(s))\,,
\end{eqnarray*}
(cf. (\ref{7}) below). Applying Theorem~1.2, we obtain the  asymptotic
formula for $\log \M^*(s)$:

\begin{thm}
Let the majorant $\M(\xi)$ satisfy conditions (\ref{2.7}), (\ref{2.9})
and (\ref{N}). Then for $s\to 0$
$$
\log \M^*(s) = (1+o(1))Q(s)\,,
$$
where $Q(s)$ is the upper Legendre transform (\ref{1.2}).
\end{thm}

\medskip As above, under an additional assumption, the remainder can be
improved to $O(\log Q(s))$.

\section{Preliminaries}
\setcounter{equation}{0}
In this section, we collect various elementary estimates which will be
used in the proofs of Theorems~1.1 and 1.2.

\subsection{The logarithmic weights}
Here, we establish several simple facts about the
logarithmic weight functions $\phi (t)$ satisfying conditions (i)--(iv), 
namely:
\begin{equation}
\lim_{t\to\infty} \frac{\phi (t)}{t} = 0\,,
\label{4.1}
\end{equation}

\begin{equation}
\lim_{t\to\infty} \phi'(t) = 0\,,
\label{4.2}
\end{equation}

\begin{equation}
\int_1^\infty \frac{\phi'(t)}{t}\, dt < \infty\,,
\label{4.3}
\end{equation}

\begin{equation}
\int_1^\infty |\phi''(t)|\, dt < \infty\,.
\label{4.4}
\end{equation}

First, let us prove relation (\ref{4.2}), then (\ref{4.1}) follows by
integration. Assume, for example, that (iv-a) holds but (\ref{4.2}) fails,
and
set $\Phi (\tau) = \phi(e^\tau)$. Then for a sequence $\tau_j\uparrow
+\infty$, $\Phi'(\tau_j)e^{-\tau_j} \ge c >0$, and
$$
\Phi(\tau_j+\xi) = \Phi(\tau_j) + \int_{\tau_j}^{\tau_j+\xi} \Phi'(\tau)\,
d\tau \ge c\xi e^{\tau_j}\,.
$$
Assuming without loss of generality that the intervals $[\tau_j,
\tau_j+1]$ are disjoint, we arrive at the contradiction:
\begin{eqnarray*}
\infty \stackrel{(iii)}> \int_1^\infty \Phi(\tau)e^{-\tau}\, d\tau 
&\ge& \sum_{j\ge 1} \int_{\tau_j}^{\tau_j+1} \Phi(\tau)e^{-\tau}\, d\tau 
\\ \\
&\ge& c \sum_{j\ge 1} \int_0^1 \xi e^{-\xi}\, d\xi = +\infty\,.
\end{eqnarray*}

The second case, when condition (iv-b) holds, is even simpler: since the
derivative $\phi'$ decreases, $\phi'(t) \ge ct$ for all $t\ge 0$, and
therefore, $\phi (t) \ge ct$ which again contradicts condition (iii).

Now, let us prove relation (\ref{4.3}): for each $A>1$
$$
\int_1^A \frac{\phi(t)}{t^2}\, dt = \phi (1) - \frac{\phi (A)}{A}
+ \int_1^A \frac{\phi'(t)}{t}\, dt\,.
$$
Letting $A\to\infty$ and using (\ref{4.1}), we obtain (\ref{4.3}).

In order to obtain (\ref{4.4}), first observe that due to (\ref{4.2}) the
integral
\begin{equation}
\int_1^\infty \phi''(t)\, dt
\label{4.5}
\end{equation}
is convergent. If condition (iv-b) holds, then the integrand in
(\ref{4.5}) is non-positive and there is nothing to prove. 
Furthermore, convergence of the integrals (\ref{4.3}) and
(\ref{4.5}) yields convergence of
$$
\int_1^\infty \frac{t^2\phi''(t)+t\phi'(t)}{t^2}\, dt\,.
$$ 
If condition (iv-a) holds, then the last integrand is non-negative,
which yields the absolute convergence of the integral (\ref{4.4}).

\subsection{Differentiating the Poisson integral}
Let
$$
u(z) = \log|W(z)| = \frac{{\rm Im}\, z}{\pi}
\intR \frac{\phi (t)\, dt}{|t-z|^2}\,,
\qquad {\rm Im}\, z>0,
$$
and  
$$
q(y) = u(iy) = \log|W(iy)|\,.
$$
In this section we collect rather straightforward estimates of the
derivatives of these two functions which will be used later.

\begin{lemma}
Let $\phi $ be a $C^2$-function
satisfying conditions (i)--(iv). Then 
\begin{equation}
u_x'(x,y) = \frac{4xy}{\pi} \int_0^\infty
\frac{t\phi'(t)\, dt}
{\left[(t-x)^2+y^2 \right]\,\left[(t+x)^2+y^2 \right]} \,,
\label{5.2}
\end{equation}
\begin{equation}
q''(y) = - \frac{4y}{\pi} \int_0^\infty \frac{t\phi'(t)\,
dt}{(t^2+y^2)^2} 
\label{5.3a}
\end{equation}
\begin{equation}
\qquad 
= \frac{2}{\pi y} \int_0^\infty 
\frac{t^2\phi''(t)}{t^2+y^2}\, dt\,,
\label{5.3b}
\end{equation}
and
\begin{equation}
u_{\theta \theta}''(r,\theta) = 
- \frac{r\sin\theta}{\pi} \intR 
\frac{t^2\phi''(t)+t\phi'(t)}{|t-re^{i\theta}|^2}\, dt\,.
\label{5.4}
\end{equation}
\end{lemma}

\smallskip\par\noindent{\em Proof: }
For $x\ge 0$ we have
\begin{eqnarray}
u(x,y) &=& \frac{y}{\pi}\int_0^\infty \phi (t)
\left[\frac{1}{(t-x)^2+y^2} + \frac{1}{(t+x)^2+y^2}
\right] dt \nonumber \\ \nonumber \\ 
&=& -\frac{1}{\pi} \int_0^\infty \phi '(t) 
\left[ \arctan \frac{t-x}{y} + \arctan\frac{t+x}{y}
\right]dt \,.
\label{5.5}
\end{eqnarray}
Integration by parts is justified since for $t\to\infty$
$$
\arctan\frac{t-x}{y} + \arctan\frac{t+x}{y} = \arctan
\frac{2ty}{x^2+y^2-t^2}
= - \frac{2y}{t} + O\left(\frac{1}{t^2}\right)\,.
$$
Differentiating (\ref{5.5}) under the integral sign, we obtain
\begin{eqnarray*}
u'_x(x,y) &=& \frac{y}{\pi} \int_0^\infty \phi'(t)
\left[ \frac{1}{(t-x)^2+y^2} - \frac{1}{(t+x)^2+y^2}
\right]dt \\ \\
&=& \frac{4xy}{\pi} \int_0^\infty 
\frac{t\phi'(t)\, dt}{\left[(t-x)^2+y^2\right]\,
\left[(t+x)^2+y^2\right]}\,,
\end{eqnarray*}
proving (\ref{5.2}).

In order to get (\ref{5.3a}) we differentiate (\ref{5.2}) by $x$ and
recall that $u''_{yy}=-u''_{xx}$. Now, since the RHS of
(\ref{5.3a}) equals 
$$
\frac{2}{\pi y} \int_0^\infty \phi'(t) d \left( \frac{t^2}{t^2+y^2}
\right)
$$
relation (\ref{5.3b}) follows by integration by parts. 

Next, we compute $ru'_r(r,\theta)$ starting again with relation
(\ref{5.5}) rewritten in the polar coordinate and then differentiating it
with respect to $r$. We have 
$$
u(r,\theta) = - \frac{1}{\pi} \int_0^\infty \phi'(t)
\left[
\arctan\left(\frac{t}{r\sin\theta}-\cot\theta\right)
+ \arctan\left(\frac{t}{r\sin\theta}+\cot\theta\right)
\right]dt\,,
$$
and
$$
\Big[ \quad . \quad \Big]'_r = -t\sin\theta 
\left[ 
\frac{1}{|t-re^{i\theta}|^2} + \frac{1}{|t+re^{i\theta}|^2}
\right]\,;
$$
therefore
\begin{equation}
ru'_r(r,\theta) 
= \frac{r\sin\theta}{\pi}
\intR \frac{t\phi'(t)\, dt}{|t-re^{i\theta}|^2}\,.
\label{5.6}
\end{equation}
Differentiation under the integral sign is justified by (\ref{4.3}). 

Iterating this procedure once more (and using (\ref{4.4}) for
justification of the differentiation), we obtain (\ref{5.4}):
\begin{equation}
- u''_{\theta \theta} = r \left( ru_r'(r,\theta) \right)'_r
=\frac{r\sin\theta}{\pi} \intR \frac{t^2\phi ''(t) + t\phi 
'(t)}{|t-re^{i\theta}|^2}\, dt\,.
\label{5.7}
\end{equation}
This completes the proof.
\footnote{Here is another argument suggested by I.~Ostrovskii. 
To prove (\ref{5.2}), observe that
\begin{eqnarray*}
u'_x(x,y) &=&  \frac{y}{\pi}
\int_{-\infty}^{\infty}
\phi(t)\, \frac{\partial}{\partial x}\frac{1}{(x-t)^2+y^2}\,dt \\
&=& 
-\frac{y}{\pi}\int_{-\infty}^{\infty}\phi(t)\,
\frac{\partial}{\partial t}\frac{1}{(x-t)^2+y^2}\,dt \\
&=&
\frac{y}{\pi}
\int_{-\infty}^{\infty}\frac{\phi'(t)\,dt}{(x-t)^2+y^2}.
\end{eqnarray*}
Relations (\ref{5.3a}) and (\ref{5.3b}) follow at once, as above. 

To prove (\ref{5.7}), we make a change of vairable in the Poisson formula
for $u(z)$: set
$z=e^{\zeta}$, $t=e^{\tau}$ for $t>0$ and $t=-e^{\tau}$ for $t<0$.
We obtain the Poisson formula for a strip:
$$
u(e^{\zeta}) =
\frac{\sin\theta}{2\pi}
\int_{-\infty}^{\infty}
\frac{\phi(e^{\tau})\,d\tau}{\cosh(\xi-\tau) -\cos\theta}
+ \frac{\sin\theta}{2\pi}
\int_{-\infty}^{\infty}
\frac{\phi(-e^{\tau})\,d\tau}{\cosh(\xi-\tau) +\cos\theta},
$$
$$
\zeta=\xi+i\theta, \quad -\infty < \xi < \infty, \quad 0<\theta<\pi\,.
$$
Differentiate the last formula twice by $\xi$.
Obviously, the second derivatives of the integrands by $\xi$ can be
replaced by their second derivatives by $\tau$. Then integrating twice by
parts we obtain
$$
[u(e^{\zeta})]^{''}_{\zeta\zeta} = 
\frac{\sin\theta}{2\pi} \int_{-\infty}^{\infty}
\frac{[\varphi(e^{\tau})]^{''}_{\tau\tau}\,d\tau}{\cosh(\xi-\tau)
-\cos\theta} +
\frac{\sin\theta}{2\pi} \int_{-\infty}^{\infty}
\frac{[\varphi(-e^{\tau})]^{''}_{\tau\tau}\,d\tau}{\cosh(\xi-\tau)
+\cos\theta}\,.
$$
Returning to the original variables, we get (\ref{5.7}) and hence
(\ref{5.4}). $\Box$
}
$\Box$

\medskip Observe several immediate corollaries:

\begin{cor}
The function $x\mapsto u(x,y)$ increases with $x$ on
$[0,\infty)$.
\end{cor}

\begin{cor}
The function $q'(y)$ decreases (that is, $q(y)$ is concave), and 
$\lim_{y\to \infty} q'(y) = 0$.
\end{cor}

\begin{cor}
If $\phi$ satisfies conditions (i)--(iii) and (iv-a), then the function 
$\theta\mapsto u(r,\theta)$ is concave on $[0,\pi]$ and
\footnote{
A. Borichev brought our attention to
the fact that the result ceases to hold if assumption (iv-a) is
replaced by (iv-b). Indeed, let $\phi (t) = \min(t, R)$. Then
\begin{eqnarray*}
u(iR) &=& 
\frac{2R}{\pi} \int_0^R \frac{t\,dt}{t^2+R^2} + \frac{2R}{\pi}
\int_R^\infty \frac{R\, dt}{t^2+R^2} 
\\ \\
&<& \frac{2R}{\pi} 
\left\{  \frac{1}{R^2}\cdot \frac{R^2}{2} + R\cdot \frac{1}{R} \right\}
=  \frac{3R}{4\pi} < \frac{R}{4} = \frac{u(\pm R)}{4}\,.
\end{eqnarray*}
}
$$
u\left(r, \frac{\pi}{2} \right) = \max_{0\le \theta \le \pi} u(r, 
\theta)\,.
$$
\end{cor}

\begin{cor}
If $\phi $ satisfies conditions (i)--(iii) and (iv-b), then the function
$q''(y)$ increases to $0$ as $y\to\infty$.
\end{cor}

\medskip In the rest of this subsection, we shall estimate the partial 
derivatives of the function $u$.

\begin{lemma}
If the function $\phi (t)$ 
satisfies conditions (i)--(iv), then
\begin{equation}
\frac{\phi'(y)}{3\pi y} \le |q''(y)| \le \frac{24}{\pi y} 
\int_y^\infty \frac{\phi (t)}{t^2}\, dt\,.
\label{5.8}
\end{equation}
In particular, 
$$
\lim_{y\to\infty} |q''(y)| = 0\,,
$$
and
\begin{equation}
\lim_{y\to\infty} y |q''(y)|^{1/2} = \infty\,.
\label{5.10}
\end{equation}
\end{lemma}

\smallskip\par\noindent{\em Proof: } We start with the lower bound in 
(\ref{5.8}). First, assume that condition (iv-a) holds. Then 
$$
|q''(y)| \stackrel{(\ref{5.3a})}\ge \frac{4y}{\pi}
\int_y^\infty \frac{t\phi'(t)}{4t^4}\, dt 
\stackrel{(iv-a)}\ge \frac{y}{\pi}\,\cdot\, y\phi'(y) \int_y^\infty
\frac{dt}{t^4} = \frac{\phi'(y)}{3\pi y} \,.
$$
If condition (iv-b) holds, then we have even a slightly better
bound:
$$
|q''(y)| \stackrel{(\ref{5.3a})}\ge \frac{4y}{\pi}
\int_0^y \frac{t\phi'(t)}{4y^4}\, dt 
\stackrel{(iv-b)}\ge \frac{\phi'(y)}{\pi y^3} \int_0^y
t\,dt = \frac{\phi'(y)}{2\pi y} \,.
$$

Now, we prove the upper bound in (\ref{5.8}) (this does not need condition 
(iv)\,):
\begin{eqnarray*}
|q''(y)| &\le&
\frac{4y}{\pi} 
\left\{
\int_0^y \frac{t\phi'(t)}{y^4}\, dt 
+ \int_y^\infty \frac{t\phi'(t)}{t^4}\, dt
\right\} \\ \\ 
&\le& \frac{4y}{\pi}
\left\{
\frac{\phi (y)}{y^3} + 3 \int_y^\infty \frac{\phi(t)}{t^4}\, dt
\right\} \\ \\
&\le& \frac{24y}{\pi} \int_y^\infty \frac{\phi (t)}{t^4}\, dt 
\le \frac{24}{\pi y} \int_y^\infty \frac{\phi (t)}{t^2}\, dt\,,
\end{eqnarray*}
completing the proof. $\Box$

\begin{lemma}
The following estimates hold:
\begin{equation}
u'_x(x,y) \ge \frac{x}{6} |q''(y)|\,,
\qquad 0\le x \le y\,;
\label{5.11}
\end{equation}
\begin{equation}
u'_x(x,y) \ge \frac{2}{5\pi} \, \min_{t\in [x/2, 3x/2]} \phi'(t)\,,
\qquad x\ge y\,;
\label{5.12}
\end{equation}
and
\begin{equation}
u'_x(x,y) \ge \frac{10}{x}\,,
\qquad x\ge 8|q''(y)|^{-1/2}\,, \quad 
y\ge y_0\,.
\label{5.13}
\end{equation}
\end{lemma}

\smallskip\par\noindent{\em Proof:} First, let $0\le x\le y$. Then
for $t\ge 0$
$$
(t-x)^2+y^2 \le 2(t^2+y^2)\,, \qquad {\rm and} \qquad
(t+x)^2+y^2 \le 3(t^2+y^2)\,.
$$
Hence
$$
u'_x (x,y) \stackrel{(\ref{5.2})}\ge 
\frac{2xy}{3\pi} \int_0^\infty \frac{t\phi'(t)\,
dt}{(t^2+y^2)^2} 
\stackrel{(\ref{5.3a})}= \frac{x}{6} |q''(y)|\,,
$$
proving (\ref{5.11}).

Now, let $x\ge y$ and $|t-x|\le \frac{y}{2}$. Then
$$
(t-x)^2+y^2 \le \frac{5y^2}{4}\,,
\qquad {\rm and} \qquad 
(t+x)^2+y^2 \le \frac{15x^2}{2}\,.
$$
Therefore,
\begin{eqnarray*}
u'_x(x,y) &\stackrel{(\ref{5.2})}\ge& 
\frac{4xy}{\pi} 
\int_{|t-x|\le y/2} \frac{t\phi'(t)\, dt}{x^2y^2}\, \cdot \frac{8}{75} \\
\\
&\ge& \frac{1}{5\pi xy} \min_{|t-x|\le y/2} \phi'(t)
\, \left[ \left(x+\frac{y}{2}\right)^2 - \left(x-\frac{y}{2}\right)^2
\right] \\ \\
&\ge& \frac{2}{5\pi}\, \min_{x/2 \le t \le 3x/2} \phi'(t)\,,
\end{eqnarray*}
proving (\ref{5.12}).

At last, let $8|q''(y)|^{-1/2}\le x$. If $x\le y$, then 
$$
u'_x(x,y) \stackrel{(\ref{5.11})}\ge \frac{x}{6} |q''(y)| \ge
\frac{10}{x}\,,
$$
and if $x\ge y$, then 
$$
u'_x(x,y) \stackrel{(\ref{5.12})}\ge 
\frac{2}{5\pi} \min_{t\in [x/2, 3x/2]} \frac{t\phi'(t)}{t} 
\ge \frac{4}{15\pi x} \min_{t\ge x/2} t\phi'(t) \ge \frac{10}{x}\,, 
$$
provided that $y\ge y_0$, where $y_0$ is large enough. This competes the proof 
of the lemma. $\Box$ 

\subsection{The Legendre transform $Q(s)$}
Here, for the convenience of references, we collect several elementary
facts about the behaviour of the function $Q(s)=\sup_{y>0}[q(y)-sy]$. 

Let $y_s$ be a unique solution of the equation
\begin{equation}
q'(y) = s\,.
\label{6.1}
\end{equation}
Due to the Corollary~3.2, $y_s\uparrow +\infty$ when $s\to 0$.
Then $Q(s)=q(y_s)-sy_s$. Differentiating this equation, we obtain
\begin{equation}
Q'(s) = q'(y_s) \frac{dy_s}{ds} - y_s - s\frac{dy_s}{ds}
\stackrel{(\ref{6.1})}= -y_s\,.
\label{6.6}
\end{equation}
Then differentiating (\ref{6.1}) by $s$, and using (\ref{6.6}), we get
$$1 = q''(y_s) \frac{dy_s}{ds} = -q''(y_s) Q''(s)\,,$$ so that
\begin{equation}
Q''(s) = - \frac{1}{q''(y_s)}\,.
\label{6.7}
\end{equation}

Now, 
\begin{equation}
Q(s) = q(y_s) - y_sq'(y_s) = \int_0^{y_s} \xi |q''(\xi)|\, d\xi
\stackrel{(\ref{5.8})}\ge \frac{\phi (y_s)}{3\pi}\,. 
\label{6}  
\end{equation}
Next, for $s\le s_0$, 
\begin{equation}
0\le \log|Q'(s)| = \log y_s \stackrel{(ii)}= o(\phi (y_s) ) 
\stackrel{(\ref{6})}= o(Q(s))\,,
\label{7}
\end{equation}
and
$$
Q''(s) = \frac{1}{|q''(y_s)|} \stackrel{(\ref{5.8})}\le
\frac{3\pi y_s}{\phi'(y_s)} \stackrel{(iv)}\le Cy_s^2\,.
$$
Hence,
\begin{equation}
0 \le \log Q''(s) = O(\log y_s) = o(\phi (y_s)\,) = o(Q(s))\,,
\qquad s\to 0\,.
\label{8}
\end{equation}

If $\phi $ has a positive lower order (condition (vi)\,), then
the estimates are much better: 
$$
\log y_s = O(\log \phi (y_s)\, ) = O(\log Q(s)\,)\,,
$$
whence
\begin{equation}
0\le \log |Q'(s)| \le O(\log Q(s)\,)\,
\label{7bis}
\end{equation}
and
\begin{equation}
0\le \log Q''(s) \le O(\log Q(s)\,)\,
\label{8bis}
\end{equation}
for $s\to 0$.

\section{Proof of Theorem 1.1}
\setcounter{equation}{0}

\subsection{The upper bound for $\left|(\F^{-1} f)(s)\right|$}
 
\begin{lemma}
Let the logarithmic weight satisfy conditions (i)-(iv). Then for $s\le
s_0$
$$
\rpW (s) \le C\, \frac{\sqrt {Q''(s)}\, e^{-Q(s)}}{|Q'(s)|^{1/p}}\,.
$$
\end{lemma}

\smallskip\par\noindent{\em Proof: } 
Let $f\in H^p(W)$, and $||f||_{H^p(W)}\le 1$.
Then by a well-known estimate
\begin{equation}
|f(z)| \le \frac{|W(z)|^{-1}}{(\pi {\rm Im} z)^{1/p}}\,,
\qquad z\in \C_+
\label{3.1}
\end{equation}
(see e.g. \cite{Koosis2})\,.
We estimate the Fourier transform $(\F^{-1}f)(s)$. 

By Cauchy's theorem, for each $y>0$ and $s>0$, 
$$
\F^{-1}f (s) = \frac{1}{\sqrt{2\pi}} \intR e^{-is(x+iy)} f(x+iy)\, dx\,.
$$
Making use of the bound (\ref{3.1}) and, 
as above, denoting by $u=\log|W|$ the
Poisson integral of the function $\phi$, $q(y)=u(0,y)$, we have
\begin{eqnarray*}
| (\F^{-1}f)(s) | &\le& 
\frac{e^{sy}}{\sqrt{2\pi}\, (\pi y)^{1/p}}
\intR e^{-u(x,y)}\, dx \\ \\
&=& \frac{e^{sy-q(y)}}{\sqrt{2\pi}\, (\pi y)^{1/p}}
\intR e^{-[u(x,y)-u(0,y)]}\, dx\,. 
\end{eqnarray*}

We choose here $y=y_s$, where $y_s$ is a unique solution of the equation
(\ref{6.1}). Then according to the definition (1.2) of the upper Legendre
transform $Q$
\begin{equation}
|(\F^{-1}f)(s)| \le  \frac{e^{-Q(s)}}{\sqrt{2\pi}\, (\pi y_s)^{1/p}}
\intR e^{-[u(x,y_s)-u(0,y_s)]}\, dx\,. 
\label{6.2}
\end{equation}
For $x\le y$,
$$
u(x,y)-u(0,y) = \int_0^x u'_x(\xi, y) \, d\xi 
\stackrel{(\ref{5.11})}\ge
\frac{1}{6} |q''(y)|\, \int_0^x \xi\, d\xi = \frac{|q''(y)|}{12} x^2\,,
$$
and therefore
\begin{eqnarray}
\int_{|x|\le y_s} e^{-[u(x,y_s)-u(0,y_s)]}\, dx
&\le& \intR \exp \left[ -\frac{1}{12} |q''(y_s)| x^2\right]\, dx
\nonumber \\ \nonumber \\
&=&\sqrt{\frac{12\pi}{|q''(y_s)|}}\,.
\label{6.3}
\end{eqnarray}

Next, for $x\ge y$, $y\ge y_0$, we have 
$8|q''(y)|^{-1/2} \stackrel{(\ref{5.10})}\le y$, and
$$
u(x,y) - u(0,y) \ge \int_{8|q''(y)|^{-1/2}}^x u'_x(\xi, y)\, d\xi
\stackrel{(\ref{5.13})}\ge 10 \log \frac{x}{8|q''(y)|^{-1/2}}\,.
$$
Therefore,
\begin{eqnarray}
\int_{|x|\ge y_s}
e^{-[u(x,y_s)-u(0,y_s)]}\, dx &\le& 2 \int_{y_s}^\infty 
\left[ \frac{8|q''(y_s)|^{-1/2}}{x}\right]^{10} dx 
\nonumber \\ \nonumber \\
&=& \frac{2\, \cdot \, 8^{10}}{9} \,
\frac{|q''(y_s)|^{-1/2}}{\left[ y_s |q''(y_s)|^{1/2}\right]^9}
\nonumber \\ \nonumber \\
&\stackrel{(\ref{5.10})}\le& |q''(y_s)|^{-1/2}\,,
\label{6.4}
\end{eqnarray}
if $s\le s_0$.

Combining estimates (\ref{6.2})--(\ref{6.4}), we obtain
\begin{equation}
|(\F^{-1}f)(s)| \le \frac{C e^{-Q(s)}}{y_s^{1/p} |q''(y_s)|^{1/2}}\,.
\label{6.5}
\end{equation}
It remains to plug in relations (\ref{6.6}) and (\ref{6.7}) in
(\ref{6.5}). $\Box$

\subsection{Asymptotics of $\left( \F^{-1}\frac{1}{W}\right)(s)$}

\begin{lemma} Let the logarithmic weight $\phi$ satisfy conditions
(i)--(v), and let 
$$
f(z) = \frac{1}{(1-iz)^{2/p} W(z)}\,,
\qquad 1\le p \le \infty\,,
$$
where the branch of the function $(1-iz)^{2/p}$ is positive when
$z=iy$, $y>0$. Then 
$$
\left( \F^{-1}f \right)(s) = (1+o(1))\, 
\frac{\sqrt{Q''(s)}\, e^{-Q(s)}}{Q'(s)^{2/p}}\,.
$$
\end{lemma}

In particular, if $p=\infty$, we get the asymptotic relation
(\ref{**}).

\smallskip\par\noindent{\em Proof: }
Set
$$
h(z) = \log W(z) = \frac{1}{\pi i} 
\intR \left[ \frac{1}{t-z} - \frac{t}{t^2+1}\right]
\phi (t)\, dt\,.
$$
Then, applying Cauchy's theorem, we get
\begin{eqnarray*}
(\F^{-1}f)(s) &=& \frac{1}{\sqrt{2\pi}} 
\intR \frac{e^{-is(x+iy) - h(x+iy)}}{(1+y-ix)^{2/p}}\, dx 
\\ \\
&=& \frac{e^{sy-q(y)}}{\sqrt{2\pi}} 
\intR \frac{e^{-isx - [h(x+iy)-h(iy)]}}{(1+y-ix)^{2/p} }\, dx
\end{eqnarray*}
since $h(iy)=q(y)$. As in the previous section, we choose here $y=y_s$
(see (\ref{6.1})\,) and split the integral into the main part
$$
I(y) \stackrel{def}= 
\int_{|x|\le \omega(y) |q''(y)|^{-1/2}}\,
\frac{e^{-isx - [h(x+iy)-h(iy)]}}{(1+y-ix)^{2/p} }\, dx
$$
and the tail
$$
J(y) \stackrel{def}= 
\int_{|x|> \omega(y) |q''(y)|^{-1/2}}\,
\frac{e^{-isx - [h(x+iy)-h(iy)]}}{(1+y-ix)^{2/p} }\, dx\,,
$$
where 
\begin{equation}
\lim_{y\to\infty} \omega (y) = \infty\,.
\label{7.2}
\end{equation}
Later, we impose other restrictions on the function $\omega (y)$.

First, we estimate the tail $J(y)$. We have 
$$
|J(y)| \le 2 \int_{\omega (y) |q''(y)|^{-1/2}}^\infty
\frac{e^{-[u(x,y)-u(0,y)]}}{(1+y)^{2/p}}\, dx\,. 
$$
If $s\le s_0$ and $y=y_s$, then $\omega (y) \ge 8$, so that in the
latter integral $x\ge 8|q''(y)|^{-1/2}$, and we can use estimate 
(\ref{5.13}):
$$
u(x,y) - u(0,y) \ge 
\int_{8|q''(y)|^{-1/2}}^x \frac{10}{\xi}\, d\xi 
= \log \left[ \frac{x}{8|q''(y)|^{-1/2}} \right]^{10}\,.
$$
Then
\begin{eqnarray*}
|J(y)| &\le& 2 \int_{\omega (y)|q''(y)|^{-1/2}}^\infty
\left[ \frac{8|q''(y)|^{-1/2}}{x} \right]^{10} 
\frac{dx}{(1+y)^{2/p}} 
\\ \\ 
&\le& \frac{2 \cdot 8^{10}}{9}\, \frac{|q''(y)|^{-5}}{\left[ 
\omega (y) |q''(y)|^{-1/2} \right]^9}\, 
\frac{1}{(1+y)^{2/p}} 
\\ \\
&\stackrel{(\ref{7.2})}=& \frac{o(1)}{|q''(y)|^{1/2}} \,  
\frac{1}{(1+y)^{2/p}}\,,  
\qquad y\to \infty\,. 
\end{eqnarray*}

Now, we consider the main integral 
$$
I(y) = \int_{|x|\le \omega (y) |q''(y)|^{-1/2}}\frac{
\exp\left\{- isx - h'(iy)x - h''(iy) \frac{x^2}{2}  
- \sum_{k=3}^\infty \frac{h^{(k)}(iy)}{k!} x^k  \right\}}
{(1+y-ix)^{2/p}}\,dx\,.
$$
Because of the choice of $y=y_s$,
$$
-is - h'(iy) = -is + iu'_y(0,y) = -i (s-q'(y)) = 0\,.
$$
So taking into account that $h''(iy)=-u''_{yy}(0,y) = -q''(y)$, we obtain
\begin{eqnarray}
I(y) &=& 
\int_{|x|\le \omega (y)|q''(y)|^{-1/2}}
\frac{\exp\left\{-|q''(y)| \frac{x^2}{2}
- \sum_{k=3}^\infty \frac{h^{(k)}(iy)}{k!} x^k  
\right\}}{(1+y-ix)^{2/p}}\, dx 
\nonumber \\ \nonumber \\
&=& \frac{1+o(1)}{y^{2/p} |q''(y)|^{1/2}}
\int_{|\xi| \le \omega (y)}
\exp \Big\{-\frac{\xi^2}{2} 
\nonumber \\ \nonumber \\ 
&\ & \qquad \qquad \qquad - \sum_{k=3}^\infty \frac{h^{(k)}(iy)}{k!} 
\frac{\xi^k}{|q'' (y)|^{k/2}} \Big\} 
d\xi\,,  
\label{7.4}
\end{eqnarray}
assuming that 
\begin{equation}
\omega (y) |q''(y)|^{-1/2} = o(y)\,,
\qquad y\to\infty
\label{7.5}
\end{equation}
(due to (\ref{5.10}) this assumption does not contradict to
(\ref{7.2})\,). 

Now, we shall show that
\begin{equation}
\sup_{|\xi|\le \omega (y)}  
\left| \sum_{k\ge 3} 
\frac{h^{(k)}(iy)}{k!} \, \frac{\xi^k}{|q''(y)|^{k/2}}\right| 
= o(1)\,, \qquad y\to\infty\,.
\label{7.7}
\end{equation}
For this, we need estimates of $h^{(k)}(iy)/k!$. We have
\begin{eqnarray*}
\left| \frac{h^{(k)}(iy)}{k!} \right| 
&=& \left| \frac{1}{\pi} \intR \frac{\phi (t)\, dt}{(t-iy)^{k+1}} \right|
\\  \\
&\le& \frac{2}{\pi} \int_0^\infty \frac{\phi (t)\,
dt}{(t^2+y^2)^{(k+1)/2}}
\\ \\
&\le & \frac{2}{\pi} 
\left\{ \frac{1}{y^{k+1}} \int_0^y \phi (t)\, dt 
+ \int_y^\infty \frac{\phi (t)}{t^{k+1}}\, dt
\right\}
\end{eqnarray*}
Then the LHS of (\ref{7.7}) is
\begin{eqnarray}
&\le& \frac{2}{\pi} 
\left\{  \frac{1}{y} \int_0^y \phi (t)\, dt\, 
\sum_{k\ge 3} \left(\frac{\omega (y)}{y|q''(y)|^{1/2}}\right)^k
+ \int_y^\infty \frac{\phi (t)}{t} 
\sum_{k\ge 3} \left(\frac{\omega (y)}{t|q''(y)|^{1/2}}\right)^k \,dt
\right\} 
\nonumber \\ \nonumber \\
&\stackrel{(\ref{7.5})}\le&
\frac{4}{\pi} \, 
\left( \frac{\omega (y)}{|q''(y)|^{1/2}} \right)^3
\left\{ \frac{1}{y^4} \int_0^y \phi (t)\, dt
+ \int_y^\infty \frac{\phi (t)}{t^4}\, dt
\right\}
\nonumber \\ \nonumber \\
&\le& \frac{16}{3\pi}  
\left( \frac{\omega (y)}{|q''(y)|^{1/2}} \right)^3
\int_y^\infty \frac{\phi (t)}{t^4}\, dt 
\nonumber \\ \nonumber \\
&\stackrel{(\ref{5.8})}\le& 
{\rm const}\, \omega^3(y) 
\frac{y^3 \int_y^\infty \frac{\phi (t)}{t^4}\,dt}{\left(
y\phi'(y) \right)^{3/2}}\,. 
\label{temp}
\end{eqnarray}
In the next to the last inequality, we used the estimate
$$
\frac{1}{y^4} \int_0^y \phi (t)\, dt 
\le \frac{\phi (y)}{y^3} 
\le \frac{1}{3} \int_y^\infty \frac{\phi (y)}{t^4}\, dt
\le \frac{1}{3} \int_y^\infty \frac{\phi (t)}{t^4}\, dt\,.
$$
By virtue of condition (v), we can choose $\omega (y)$ increasing to
infinity with $y$ and such that the RHS of (\ref{temp}) is $o(1)$.
This proves (\ref{7.7}).

Now, making use of (\ref{7.7}), we continue the estimate
(\ref{7.4}):
$$
I(y) = \frac{1+o(1)} {y^{2/p} |q''(y)|^{1/2}} 
\int_{|\xi|\le \omega (y)} e^{-\xi^2/2}\, d\xi 
= \frac{\sqrt{2\pi} + o(1)}{y^{2/p} |q''(y)|^{1/2}}\,,
$$
completing the proof. $\Box$

\subsection{Conclusion of the proof of Theorem~1.1}
Juxtaposing lemmas~4.1 and 4.2, we get
\begin{eqnarray*}
&\,& \qquad -Q(s) + \frac{1}{2}\log Q''(s) - \frac{2}{p}\log|Q'(s)| - O(1)
\\ \\ &\,& \qquad \qquad 
\le \log \rpW (s) 
\le -Q(s) + \frac{1}{2}\log Q''(s) - \frac{1}{p}\log|Q'(s)| + O(1)\,.
\end{eqnarray*}
Now, estimates (\ref{7}) and (\ref{8}) give us (\ref{*}). 
When condition (vi) holds, 
estimates (\ref{7bis}) and (\ref{8bis}) yield  (\ref{***}). 
$\Box$

\section{Proof of Theorem~1.2}
\setcounter{equation}{0}
As above, we prove separately an upper and a lower bound for
$\lambda_W(s)$ which together immediately yield the theorem.

\subsection{A version of the Laplace asymptotic estimate}
Here, we give a version of the Laplace asymptotic estimate needed for the
proof of Theorem~1.2. This result is a minor
modification of the well-known (cf. \cite[\S 17]{R}),
however, for the reader's convenience we bring it with the proof.

\begin{thm}
Let
\begin{equation}
N(s) = \int_1^\infty y^a e^{-sy +q(y)}\, dy\,,
\label{L}
\end{equation}
where $a\in \R$, 
and $q$ is a $C^3(0,\infty)$ concave function such that 
\begin{itemize}
\item[(a)] $q$ steadily increases to $+\infty$ with $y$;
\item[(b)] $q'$ steadily decreases to $0$ as $y$ tends to $+\infty$;
\item[(c)] $q''$ steadily increases to $0$ as $y$ tends to $+\infty$;
\item[(d)] 
$$ 
\lim_{y\to\infty} y^2|q''(y)| = +\infty\,; 
$$
\item[(e)] 
$$ 
\lim_{y\to \infty} \frac{|q''(y)|^{3/2}}{q'''(y)} = +\infty\,.
$$
\end{itemize}
Then for $s\to 0$
\begin{equation}
N(s) = (1+o(1)) |Q'(s)|^a
\sqrt{2\pi Q''(s)} e^{Q(s)}\,,
\label{8.3}
\end{equation}
where 
$$
Q(s) = \sup_{y>0} [q(y) - sy]
$$
is the upper Legendre transform of $q$.
\end{thm}

\smallskip The proof is based on two auxiliary lemmas. 

\begin{lemma} In assumptions of Theorem~5.1, 
there is a function $\eta (y) \le y/2$, $1\le y<\infty$, such that 
\begin{equation}
\lim_{y\to \infty} \eta (y) \sqrt{|q''(y)|} = \infty\,,
\label{8.4}
\end{equation}
\begin{equation}
\lim_{y\to \infty} \frac{\eta (y)}{y} = 0\,,
\label{8.5}
\end{equation}
and
\begin{equation}
\lim_{t\to\infty} \, \sup_{t-\eta (t) \le y \le t+ \eta (t)}
\left| \frac{q''(y)}{q''(t)} - 1 \right| = 0\,.
\label{8.6}
\end{equation}
\end{lemma}

\smallskip\par\noindent{\em Proof of Lemma~5.1: } Let $y>t$. Then for some
$\xi\in (t,y)$
\begin{eqnarray*}
\left| \frac{q''(y)}{q''(t)} - 1 \right| &=&
(y-t) \frac{q'''(\xi)}{|q''(t)|} \\ \\
&\le& (y-t) \sqrt{|q''(t)|} \frac{q'''(\xi)}{|q''(\xi)|^{3/2}} \\ \\
&\le& (y-t) \frac{\sqrt{|q''(t)|}}{\gamma (t)}\,,
\end{eqnarray*}
where 
$$
\gamma (t) \stackrel{def}= \inf_{\xi\ge t}
\frac{|q''(\xi)|^{3/2}}{q'''(\xi)}\,,
\qquad \gamma (t)\uparrow\infty\,.
$$ 

Let $\gamma_1(t)$ be a minorant of $\gamma (t)$ such that
$\gamma_1(t)\uparrow \infty$, and $\gamma_1(t)=o(1)t^2|q''(t)|$, 
for $t\to\infty $. Defining
$$
\eta_1 (t) \stackrel{def}= \min\left(
  \sqrt{\frac{\gamma_1(t)}{|q''(t)|}}, 
\frac{t}{4} \right),
$$
we obtain a one-sided version of (\ref{8.6}):
$$
\sup_{t\le y \le t+\eta_1(t)}\, \left| \frac{q''(y)}{q''(t)} - 1 \right|
\le \eta_1 (t) \frac{\sqrt{|q''(t)|}}{\gamma(t)} 
\le \frac{\sqrt{\gamma_1(t)}}{\gamma (t)} 
\le \frac{1}{\sqrt{\gamma (t)}} \to 0\,,
$$
together with 
$$
\frac{\eta_1(t)}{t} \le \sqrt{\frac{\gamma_1(t)}{t^2|q''(t)|}}
= o(1)\,, 
$$
and
$$
\eta_1(t) \sqrt{|q''(t)|}
= \min\left( \sqrt{\gamma_1(t)}, \frac{1}{4}t|q''(t)|\right)
\to \infty\,,
$$
when $t\to\infty$.

Set $\eta (t) = \eta_1 (\tau_t)$, where
$\tau_t\stackrel{def}= \inf\{\xi:\, \xi+\eta_1(\xi)\ge t \}$.
Then $\frac{2}{3} t \le \tau_t \le t$, and (\ref{8.5}) follows. Since the
function $|q''(t)|$ decreases, relation
(\ref{8.4}) follows as well. At last, $\tau_t + \eta_1 (\tau_t) = t$, so
that
$$
\sup_{t-\eta_1(\tau_t) \le y \le t} \,
\left| \frac{q''(y)}{q''(t)} - 1\right| = 
\frac{|q''(\tau_t)|}{|q''(t)|} - 1 = 
= \frac{|q''(\tau_t)|}{|q''(\tau_t + \eta_1(\tau_r)\,)|} -1 = 
o(1), \qquad t\to\infty
$$
by the choice of the function $\eta_1(t)$. $\Box$

\smallskip The next lemma is useful for estimation of the tails.

\begin{lemma}
Let $\psi\in C^2(a,\infty)$ be an increasing convex function. Then
$$
\int_a^\infty e^{-\psi (x)}\, dx < \frac{e^{-\psi (a)}}{\psi'(a)}\,.
$$
\end{lemma}

\smallskip\par\noindent{\em Proof of the Lemma~5.2: } 
Since $\psi $ is convex and increasing, 
$$
\psi (x) - \psi (a) \ge (x-a) \psi'(a)\,,
$$
so that 
$$
\int_a^\infty e^{-\psi (x)}\, dx \le
\int_a^\infty e^{-\psi (a) - (x-a)\psi'(a)}\, dx 
= \frac{e^{-\psi (a)}}{\psi' (a)}\,.
$$
Done! $\Box$

\medskip\par\noindent{\em Proof of the Theorem~5.1 } Let $t=t_s$ be a
unique solution of the equation
$$
q'(t) = s\,,
\qquad s\le s_0 = q'(+0)\,.
$$
Then we split the integral $N(s)$ into three parts:
$$
N(s) = \left( \int_1^{t-\eta} + \int_{t-\eta}^{t+\eta} 
+ \int_{t+\eta}^\infty  \right)
y^a e^{-sy + q(y)}\, dy = \sum_{k=1}^3 I_k(s)\,,
$$
where $t=t_s$, $\eta = \eta (t_s)$ is a function from Lemma~5.1.
The asymptotics (\ref{8.3}) is
defined by the integral $I_2(s)$; the other two integrals give the
remainder. 

Let us start with the principal term:
\begin{eqnarray}
I_2(s) 
&=& e^{q(t)-st} \int_{t-\eta}^{t+\eta} y^a 
e^{q(y)-q(t)-(y-t)q'(t)}\, dy 
\nonumber \\ \nonumber \\
&=& (1+o(1)) t^a e^{q(t)-st} \int_{t-\eta}^{t+\eta}  
\exp \left\{ -\frac{1}{2} (y-t)^2 |q''(\xi (y))|\right\} \, dy 
\nonumber \\ \nonumber \\
&\stackrel{(\ref{8.6})}=& (1+o(1)) |Q'(s)|^a e^{Q(s)}
\int_{t-\eta}^{t+\eta}  
\exp \left\{ -\frac{1}{2} (y-t)^2 |q''(t)|\right\} \, dy 
\nonumber \\ \nonumber \\
&=& (1+o(1)) |Q'(s)|^a e^{Q(s)} 
\sqrt{2\pi Q''(s)}
\label{8.7}
\end{eqnarray}
for $s\to 0$.

For estimates of the integrals $I_1$ and $I_3$ we use Lemma~5.2. First, we
observe that the function $y\mapsto sy-q(y) - a \log y$ is convex. 
Therefore, the Lemma~5.2 is applicable:
\begin{eqnarray}
\int_{t+\eta}^\infty y^a e^{-sy+q(y)}\, dy 
&\le& (t+\eta)^a \frac{e^{-s(t+\eta)+ q(t+\eta)}}{s-q'(t+\eta) -
a(t+\eta)^{-1}} 
\nonumber \\ \nonumber \\
&=& (1+o(1)) t^a e^{q(t)-st} \frac{e^{q(t+\eta) - q(t) - \eta q'(t)}}
{q'(t) - q'(t+\eta) + O(1/t)} 
\nonumber \\ \nonumber \\
&=& (1+o(1)) |Q'(s)|^a e^{Q(s)} 
\frac{e^{-\frac{1+o(1)}{2} \eta^2 |q''(t)|}}
{\eta |q''(t)| (1+o(1)) + O(1/t)} 
\nonumber \\ \nonumber \\
&\le& (1+o(1)) |Q'(s)|^a \sqrt{Q''(s)} e^{Q(s)} 
\frac{e^{-\frac{1}{4} \eta^2 |q''(t)|}}
{\eta |q''(t)|^{1/2} + O(\frac{1}{t|q''(t)|^{1/2}})} 
\nonumber \\ \nonumber \\
&=& o(1) |Q'(s)|^a \sqrt{Q''(s)} e^{Q(s)} 
\qquad s\to 0\,.
\label{8.8}
\end{eqnarray}
The same estimate holds for the second integral:
\begin{equation}
I_2(s) < (t-\eta)^a \frac{e^{-s(t-\eta) + q(t-\eta)}}{q'(t-\eta) - s +
O(1/t)} =
o(1) |Q'(s)|^a \sqrt{Q''(s)} e^{Q(s)}\,. 
\label{8.9}
\end{equation}
Collecting estimates (\ref{8.7})--(\ref{8.9}), we finally obtain
(\ref{8.3}). 
$\Box$

\subsection{The upper bound for $\lW (s)$}
\begin{lemma} Let $\phi (t)$ satisfy conditions (i)--(iii), (iv-a), and 
(iv-b). Then  for $s\to 0$
\begin{equation}
\lW (s) \le (1+o(1)) \sqrt{2\pi Q''(s)}\, e^{Q(s)}\,.
\label{upper}
\end{equation}
\end{lemma}

\smallskip\par\noindent{\em Proof: }
Let $||f||_{B(W)}\le 1$. Then
$$
|f(z)|, |f({\overline z})| \le |W(z)|\,,
\qquad z\in \C_+\,.
$$
Rotating the integration line in the Laplace transform, and then using the
Corollary~3.3, we get
\begin{eqnarray*}
\left| (\L f)(se^{i \psi}) \right| &=& 
\left| \int_0^\infty f(re^{-i\psi}) e^{-sr}\, dr
\right| 
\\ \\
&\le& \int_0^\infty \exp\left\{ 
\log|W(re^{-i\psi})| - sr
\right\} \, dr
\\ \\
&\le& 
\int_0^\infty \exp\left\{ \log|W(ir)| - sr \right\} \, dr 
\\ \\
&\le& 
\int_1^\infty \exp\left\{ \log|W(ir)| - sr \right\} \, dr 
+ O(1)\,, \qquad s\to 0\,. 
\end{eqnarray*}

Now, we check that the function $q(y)=\log|W(iy)|$ meets conditions
(a)--(e) of Theorem~5.1, and then apply this theorem with
$a=0$. Conditions (a) and (b) follow from 
Corollary~3.2, (c) follows from Corollary~3.4, and condition
(d) follows from the lower bound in the estimate (\ref{5.8}) and
(iv-b). In order to check condition (e), we estimate 
from above the third
derivative $q'''$ using relation (\ref{5.3b}). 
Differentiating (\ref{5.3b}) once 
(this is permitted due to (\ref{4.4})\,), we get
\begin{eqnarray*}
q'''(y)&=& \frac{2}{\pi y^2} \int_0^\infty
\frac{t^2[-\phi''(t)]}{t^2+y^2}\, dt
+ \frac{4}{\pi} \int_0^\infty \frac{t^2[-\phi''(t)]}{(t^2+y^2)^2}\, dt \\
\\
&\le& \frac{6}{\pi y^2} \int_0^\infty \frac{t^2[-\phi''(t)]}{t^2+y^2}\, dt 
= \frac{3|q''(y)|}{y}\,,
\end{eqnarray*}
and therefore
$$
\frac{q'''(y)}{|q''(y)|^{3/2}} \le \frac{3}{y|q''(y)|^{1/2}} 
\stackrel{(\ref{5.10})}= o(1)\,,
\qquad y\to\infty\,.
$$

Thus, applying the Theorem~5.1, we get
$$
\left| (\L F)(se^{i\psi})\right|
\le (1+o(1)) \sqrt{2\pi Q''(s)}\,
e^{Q(s)}\,,
\qquad s\to 0\,.
$$
This proves the upper bound (\ref{upper}). $\Box$

\subsection{The lower bound for $\lW (s)$}
\begin{lemma} Let $\phi (t) $ satisfy conditions (i)--(iii), (iv-a), and
(iv-b). 
Then there is an entire function $E\in B(W)$ such that for some
$N<\infty$
\begin{equation}
\left|\left( \L E\right)(-is)\right|
\ge (1+o(1))\, |Q'(s)|^{-N}\, \sqrt{Q''(s)} e^{Q(s)}\,,
\qquad s\to 0\,.
\label{lower}
\end{equation}
\end{lemma}

First, applying a result of Y.~Domar \cite[Lemma~4]{Domar} (more
precisely, we use his intermediate estimate (19)), we find an even entire
function
$$
G(z) = \sum_{n\ge 0} a_{2n} z^{2n}\,,
\qquad a_{2n}\ge 0\,,
$$
such that
$$
O < c_1 |x|^{-2N} \le G(x)e^{-2\phi (x)} \le c_2 < \infty\,,
\qquad |x|\ge 1\,.
$$ 
The constants $c_1$ and $c_2$ depend on the function $\phi $ and are
independent of $x$. Since
$$
M(r, G) \left( \stackrel{def}= \max_{|z|=r} |G(z)| \right) = G(r)\,,
$$
the function $G$ has zero exponential type, and moreover, belongs to the
convergence class:
$$
\int^\infty \frac{\log M(r, G)}{r^2}\, dr < \infty\,.
$$
Applying the Krein-Akhiezer factorization theorem
\cite[Appendix~V]{Levin}, we factorize 
$$
G(z) = E(z)\, \overline{E({\overline z})}\,,
$$
where $E$ is an entire function of zero exponential type with zeroes in
the lower half-plane. Since $G(x)=|E(x)|^2$, 
\begin{equation}
c_1 |x|^{-N} \le |E(x)|e^{-\phi (x)} \le c_2\,,
\qquad |x|\ge 1\,.
\label{9.1}
\end{equation}
In particular, $E/W$ is bounded on the real axis. Since $E$ has zero
exponential type and $W$ is outer, we can apply the Phragm\'en-Lindel\"of
principle to the functions $E/W$, $E^*/W$, and conlude that these
functions are in $H^\infty$. Hence $E\in B(W)$.

Without loss of generality, we assume that $E(0)=G(0)=1$, so that $E$ is a
canonical product of
genus zero with the zero set symmetric with respect to the imaginary axis
(since $G$ is even). Thus $E(iy)\ge 0$ for $y\ge 0$.
Applying again the Phragm\'en-Lindel\"of principle to the function
$E/W$ in $\C_+$ and using estimate (\ref{9.1}), we get
\begin{equation}
E(iy) \ge c_1 y^{-N} |W(iy)|\,,
\qquad y\ge 1\,. 
\label{9.2}
\end{equation}
Therefore, for $s\to 0$, 
$$
\left| \left( \L E\right)(-is)\right| = 
\int_0^\infty E(iy) e^{-sy}\, dy
\stackrel{(\ref{9.2})}\ge c\, \int_1^\infty 
y^{-N} e^{\log|W(iy)|-sy}\, dy - O(1)\,.
$$
Above, in the proof of the previous lemma, we already checked that
the function $q(y)=\log |W(iy)|$ satisfies all assumptions of 
Theorem~5.1. Applying this result, we obtain the estimate
(\ref{lower}). $\Box$

\bigskip\par\noindent{\bf Remark } If the logarithmic weight $\phi $
satisfies an additional assumption

\smallskip\par\noindent (vii)
$$
\liminf_{\tau \to \infty} \frac{d^2\phi (e^\tau)}{d\tau^2} > 0\,,
$$
then instead of the result of Domar we may use a
recent result of Borichev
\cite{Borichev3} and construct an entire function $G$ of genus zero
with non-negative Taylor coefficients and such that
$$
0<c_1 \le G(x)e^{-2\phi (x)} \le c_2 <\infty\,,
\qquad x\in \R\,.
$$
Repeating verbatim the rest of the argument, we obtain an entire function
$E\in B(W)$ such that 
$$
\left| (\L E)(-is) \right| \le (1+o(1)) \sqrt{Q''(s)}\, e^{Q(s)}\,,
$$
whence
$$
0<a\le \frac{\lW(s)}{\sqrt{Q''(s)}\, e^{Q(s)}} \le b<\infty\,. 
$$

\bigskip\par\noindent Vladimir Matsaev: {\em School of Mathematical
Sciences,
Tel-Aviv University, \newline
\noindent Ramat-Aviv, 69978, Israel

\par\noindent matsaev@post.tau.ac.il}

\bigskip\par\noindent Mikhail Sodin: {\em School of Mathematical
Sciences, Tel-Aviv University, \newline
\noindent Ramat-Aviv, 69978, Israel

\par\noindent sodin@post.tau.ac.il}


\begin{thebibliography}{2}


\bibitem{Bang} {\sc Th. Bang},
\newblock{\em The theory of metric spaces applied to infinitely
differentiable functions,}
\newblock{Math. Scand. {\bf 1} (1953), 137--152. }

\bibitem{Beurling} {\sc A. Beurling},
\newblock{\em Extremal Distance and Estimates for Harmonic 
Measure. Mittag-Leffler Lectures on Complex
Analysis (1977-1978),}
\newblock{Collected Works, Vol \rm {I}, Birkh\"auser, Boston, 1989.}

\bibitem{Borichev} {\sc A. Borichev},
\newblock{\em Analytic quasi-analyticity and asymptotically
holomorphic functions,}
\newblock{Algebra i Analiz {\bf 4} (1992), no. 2, 70--87. English  
transl. in St. Petersburg Math. J. {\bf 4} (1993), 
259--272.}

\bibitem{Borichev2} {\sc A. Borichev},
\newblock{\em Beurling algebrals and the generalized Fourier transform,}
\newblock{Proc. London Math. Soc., {\bf 73} (1996), 431--480.}

\bibitem{Borichev3} {\sc A. Borichev},
\newblock{\em The polynomial approximation in Fock-type spaces,}
\newblock{Math. Scand. {\bf 82} (1998), 256--264.}

\bibitem{Borichev4} {\sc A. Borichev},
\newblock{\em On the closure of polynomials in weighted spaces of
functions on the real line,}
\newblock{Indiana Univ. Math. Journ. {\bf 50} (2001)}

\bibitem{Carleman} {\sc T. Carleman},
\newblock{\em Fonctions quasianalytiques,}
\newblock{Paris, Gauthier Villars, 1926.}

\bibitem{Carleman2} {\sc T. Carleman},
\newblock{\em Extension d'un th\'eor\`eme de Liouville,}
\newblock{Acta Math. {\bf 48} (1926), 363--366.}
\newblock{\'Edition Compl\'ete des Articles, Malm\"o, 1960.}

\bibitem{Domar_a} {\sc Y. Domar},
\newblock{\em On the existence of a largest subharmonic minorant of a
given function,}
\newblock{Ark. Math. {\bf 3} (1957), 429--440.}


\bibitem{Domar} {\sc Y. Domar}, 
\newblock{\em Closed primary ideals in a class of Banach algebras,}
\newblock{Math. Scand. {\bf 7} (1959), 109--125.}

\bibitem{Domar_b} {\sc Y. Domar},
\newblock{\em Uniform boundedness in families related to subharmonic
functions,}
\newblock{J. London Math. Soc. {\bf 38} (1988), 485--491.}


\bibitem{D1} {\sc E. Dyn'kin},
\newblock{\em Functions with a given estimate for $\partial
f/\overline{\partial} z$ and N.~Levinson's theorem, }
\newblock{Math. USSR Sb. {\bf 18} (1972), 181--189.}

\bibitem{D2} {\sc E. Dyn'kin},
\newblock{\em The growth of an analytic function near its
set of singular points,}
\newblock{Investigations on linear operators and the theory of functions,
III. Zap. Nau\v cn. Sem. Leningrad.
Otdel. Mat. Inst. Steklov. (LOMI) {\bf 30} (1972),
158--160. (in Russian).}

\bibitem{D3} {\sc E. Dyn'kin},
\newblock{\em The pseudoanalytic
extension,}
\newblock{J. Anal. Math. {\bf 60} (1993), 45--70.}

\bibitem{Gurarie} {\sc V. P. Gurarie}, 
\newblock{\em On Levinson's theorem concerning families of analytic
functions,}
\newblock{Investigations on linear operators and the theory of functions,
I. Zap. Nau\v cn. Sem. Leningrad.
Otdel. Mat. Inst. Steklov. (LOMI) {\bf 19} (1970),
215--220. (in Russian).}

\bibitem{HW} {\sc I. I. Hirschman and D. V. Widder},
\newblock{\em The Convolution Transform, }
\newblock{Princeton Univ. Press, Princeton, New Jersey, 1955.}

\bibitem{Koosis2} {\sc P. Koosis},
\newblock{\em Introduction to $H_p$ Spaces, 2nd ed.,}
\newblock{Cambridge Univ. Press, Cambridge, 1998.}

\bibitem{Koosis} {\sc P. Koosis},
\newblock{\em The Logarithmic Integral, {\rm I},}
\newblock{Cambridge Univ. Press, Cambridge, 1988.}

\bibitem{Levin} {\sc B. Ya. Levin},
\newblock{\em Distribution of zeros of entire functions,}
\newblock{Amer. Math. Soc., Providence, RI, 1980.}

\bibitem{Levinson} {\sc N. Levinson},
\newblock{\em Gap and density theorems,}
\newblock{Colloquium Publ. {\bf 26}. Amer. Math. Soc., 1940.}

\bibitem{Matsaev} {\sc V. Matsaev},
\newblock{\em Uniqueness, completeness and compactness theorems 
related to the classical concept of quasianalyticity.}
\newblock{PhD Thesis, Kharkov, 1964.}

\bibitem{Nik} {\sc N. Nikol'skii},
\newblock{\em Yngve Domar's Forty Years in Harmonic Analysis,}
\newblock{Acta Universitatis Upsaliensis {\bf 58} (1995), 45--78.}
 
\bibitem{R} {\sc E. Ya. Riekstyn'\v{s}},
\newblock{\em Asymptotic expansions of integrals, {\rm II},}
\newblock{Zinatne, Riga, 1977. (in Russian)}

\bibitem{Sjoberg} {\sc N. Sj\"oberg},
\newblock{\em Sur les minorantes sousharmoniques d'une fonction don\'ee,}
\newblock{Neuvi\`eme Congr\'es Math\'ematique Scandinaves 1938, 309--319.}

\bibitem{Volberg} {\sc A. Volberg},
\newblock{Seminar talk, Kharkov University function theory
seminar, Spring 1990.}

\end{thebibliography}
\end{document}